
\input amstex
\documentstyle{amsppt}

\def\t#1 {\text{\rm #1}}
\def\b#1 {\text{\rm\bf #1}}


\input BoxedEPS.tex
\SetepsfEPSFSpecial
\HideDisplacementBoxes
\def\figure#1#2{\medskip \centerline{\BoxedEPSF{fig#1.eps}}
  \smallskip\centerline{#2}\medskip}


\topmatter
\title Lectures on controlled topology: mapping cylinder neighborhoods
  \endtitle

\author  Frank Quinn \endauthor
\rightheadtext{ Controlled topology: mapping cylinders}
\address 
Mathematics, Virginia Tech, Blacksburg VA 24061-0123\endaddress
\email quinn\@math.vt.edu\endemail
\date January 2002\enddate
\thanks Partially supported by the US National Science Foundation and the ICTP\endthanks
\abstract The existence theorem for mapping
cylinder neighborhoods is discussed as a prototypical example of controlled topology and its
applications. The first of a projected series developed from lectures  at the Summer
School on High-Dimensional Topology, Trieste Italy 2001\endabstract
\endtopmatter

\head 1. Introduction\endhead
Controlled topology has the hallmarks of a mature mathematical subject: powerful results,
sophisticated interactions with, and applications to, other subjects,  difficult and
unexpectedly beautiful conjectures. It is not very accessible, however. Partly this is
because complete results are difficult and there is not a large enough community to sustain
interest in partial answers. Another problem is that it blossomed rapidly, so lacks the
more-accessible historical development and expositions of most  mature subjects. This
paper is the first in a projected series to try to address this. Here we outline 
the setting and applications of the existence theorem for mapping cylinder neighborhoods
(originally, ``completions'' of ends of maps). This illustrates most of the ingredients of the
subject: local homotopy theory, local fundamental groups, elaborate algebraic obstructions,
interesting applications. The focus is  on what all these
things mean and how they fit together, and most details are omitted. 

This paper is an expansion of the first third of a series of lectures given at the Summer
School on High-Dimensional Topology in Trieste, Italy, in the summer of 2001. Other topics
were the controlled h-cobordism theorem, illustrating some of the geometric and algebraic
techniques; and homology manifolds, illustrating the still-incomplete theory of controlled
surgery.
  
\subhead 1.1 Locating the subject\endsubhead
In the first half of the 20th century topology had two main branches: point-set topology,
concerned with local properties (separation, connectedness, dimension theory etc); and
algebraic topology, concerned with definition and detection of global structure (homology,
characteristic classes, etc.). In the 50s and 60s the algebraic branch split into homotopy
theory and geometric topology. Homotopy theory was still largely descriptive, but in the
geometric area the emphasis changed from description to construction. For instance rather than
computing homology of examples of manifolds, the objective was to construct or
classify manifolds with given homological structure. This development was mainly restricted
to spaces with uniform local structure, i.e\. manifolds. Some of the descriptive techinques
had  extensions to spaces with symmetries (group actions) and stratified spaces
such as algebraic varieties. Extensions of constructive methods were very limited due to
complexity of interactions between different levels in group actions or stratified sets.

Controlled topology began in the late 1970s and 80s as a way to apply the constructive
techniques of geometric topology to local questions more typical of point-set topology. For
example, which subspaces of a space have a neighborhood
homeomorphic to a mapping cylinder?  Mapping cylinders must be constructed
rather than simply detected, so although this is a  local  question it
requires constructive techniques. 

Controlled topology has had striking successes in elucidating geometric structure.
Unexpectedly, it has also had striking applications in algebra. Geometric problems have
obstructions related to linear or quadratic algebra ($K$- and $L$-theory). Controlled
geometric problems have obstructions in controlled-algebra analogs, essentially homology with
coefficients in spectra related to the uncontrolled obstructions.  This turns out to be a
two-way street: results about ordinary obstructions give information
about control and local geometric structure, and conversely direct controlled constructions
can give information about ordinary obstructions. 
For instance the famous ``strong Novikov conjecture'' asserts that some
 assembly maps from homology to ordinary obstruction groups are isomorphisms, at least
rationally. The homology corresponds to controlled problems, ordinary groups correspond to
uncontrolled problems, and the assembly map corresponds to ``forgetting control.'' If the
assembly map is an isomorphism then solvability of an uncontrolled problem determines
solvability of a more delicate controlled one. Conversely if there is a geometric
construction that ``gains control'' --- produces a controlled solution from an uncontrolled
one --- then the assembly map must be an isomorphism. Most geometric proofs of cases
of the Novikov conjecture rely on this principle, and the most delicate (especially work of
Farrell and Jones) use it explicitly. 

Controlled topology thus lies at the juncture between geometric and point-set topology,
homotopy theory, and stable algebra. Constructions and proofs tend to be elaborate, but
outcomes can be deep and powerful.

\subhead 1.2 Plan \endsubhead
True mastery of a subject requires understanding the details. However to get started, or for
those looking for application rather than mastery, an overview can be helpful. The goal here
is to give such an overview: definitions and enough explanation for good understanding of the
statements of theorems, sketches of proofs in enough detail to show how the hypotheses are
used and what the  difficult points are. Finally in this paper we focus on the construction
and application of mapping cylinders. This illustrates most of the techniques and issues
of controlled topology.

 The
central result in the paper is the existence theorem for mapping cylinder neighborhoods, 3.1.
However the hypotheses are quite elaborate, so Section 2 is devoted to developing them.
Specifically, 2.1 gives the definition, 2.2 describes the use of ``control spaces'', and 2.3
describes the simplest (``uncontrolled'') case of neighborhoods of points. Section 2.4 defines
tameness and describes some of the results. Tameness often does not appear in statements of
applications because it follows from other hypotheses, but it is central to proofs. Section
2.5 discusses homotopy links, which provide homotopy models and play an important role in
controlling local fundamental groups, as explained in 2.6. Stratified systems of fibrations
are introduced in 2.7. These are needed to impose some regularity on local fundamental
groups, and appear prominently in the structure of stratified sets. Section 2.8 begins
development of the ``spectral sheaf homology'' used to describe the obstruction groups. A key
feature of this theory is the assembly map defined in 2.9. The most elaborate part of the
development is in section 2.10, where the controlled $K$-theory (more precisely,
pseudoisotopy) spectrum is discussed. It is in large part the controlled assembly isomorphism
theorem for this spectrum that makes the theory accessible and useful.

With the setting and hypotheses explained, the existence theorem for mapping cylinder
neighborhoods is stated in 3.1. The proof is outlined in sections 3.2 and 3.3. Some useful
refinements are given in 3.4. The first concerns smooth and PL structures, and the second
gives a recognition criterion for mapping cylinders. 

Applications are given in section 4. The first three are straightforward: the special case of
manifolds (4.1) with its corollary the finiteness of compact finite dimensional ANRs (the
Borsuk conjecture), and collaring in homology manifolds (4.3). Then we consider mapping
cylinders in cases where the tameness and local fundamental group structure are more
elaborate, namely stratified spaces (4.4), and a special case where the obstructions can be
made relatively explicit, topological actions of finite groups (4.5). The final application
(in 4.6) is to define topological regular neighborhoods. These are mapping cylinders in a
product with $[0,\infty)$, and generalize the ``approximate tubular neighborhoods'' in
stratified spaces developed by Hughes and others.

\head 2. The setting\endhead
To illustrate the basic ideas of control we investigate the existence of mapping cylinder
neighborhoods.  Fix a space $X$ and a closed subset $Y$. We suppose the complement $X-Y$
is a manifold since we use manifold techniques there. Finally we suppose $X$ and $Y$ are
finite-dimensional ANRs (absolute neighborhood retracts), to avoid local point-set problems.

\subhead 2.1 Definition\endsubhead A {\it mapping cylinder neighborhood\/} of $Y$ in
$X$ is a closed neighborhood $N\supset Y$ with frontier $\partial N$ a submanifold of $X-Y$,
a map
$f\:\partial N\to Y$, and a homeomorphism of the mapping cylinder of $f$ with $N$ which is
the identity on $\partial N$ and $Y$. 

More explicitly the mapping cylinder is the
identification space
$\partial N
\times I \cup_f Y$, where ``$\cup_f$'' indicates that points $(x,0)\in \partial N\times I$
are identified with $f(x)\in Y$, and the homeomorphism $\partial N\times I\cup_f Y\to N$
restricts to the identities on  $\partial N\times \{1\}\to \partial N$ and $Y\to Y$. 

\subhead 2.2 The control space \endsubhead 
There is a canonical projection of a mapping cylinder to the subspace $Y$. In fact we have
assumed $Y$ is an ANR so there is  a projection of a neighborhood to $Y$ whether there is a
mapping cylinder or not. Denote this by $p\:N-Y\to Y$. We refer to $Y$ as
the control space and the projection as the  control map. To explain the terminology we
observe that a mapping cylinder is equivalent to a product structure on the complement
$N-Y\simeq \partial N\times (0,1]$ so that the images of  open arcs
$p(\{x\}\times(0,1])$ converge. The map $f\:\partial N\to Y$ can be recovered as the limit
$f(x)=\lim_{t\to 0}p(x,t)$.  Convergence is arranged using the Cauchy criterion:
constructions are done so that images of subintervals $p(\{x\}\times[\frac1{n+1},
\frac1{n}])$ have preassigned small diameter. The crucial issue is control of sizes of
images in
$Y$, hence the ``control'' terminology.

\subhead 2.3 The uncontrolled case\endsubhead 
Most controlled  theorems have older ``uncontrolled'' versions in which the control
space is implicitly taken to be a point. This version of the mapping cylinder question is:
when does a point have a neighborhood homeomorphic to a cone? Recalling that we have assumed
the complement is a manifold, this can be reformulated as: when is a noncompact manifold the
interior of a compact manifold with boundary? The 1-point compactification then plays the
role of $X$, and $Y$ is the point at infinity. Theorems of Browder,
Livesay and Levine \cite{BLL} and Siebenmann \cite{S} answer this: there is a necessary
homotopy-theoretical ``tameness'' condition, and then an obstruction in  algebraic
$K$-theory. 

We expand on the $K$-theory part. When the tameness condition is satisfied  we can
predict the fundamental group of the boundary of the neighborhood: the groups $\pi_1(U-Y)$
indexed by the inverse system of neighborhoods $U$ of $Y$ converges nicely to a finitely
presented group
$\pi$. In the course of the construction of actual boundaries a finitely presented projective
module over the group ring $Z[\pi]$ is encountered. If this projective module is stably free
the construction can be continued to give the desired structure. The obstruction is therefore
essentially the class of this module in the group of stable equivalence classes of projective
modules,
$K_0(Z[\pi])$. This is a little too big: $K_0$ records the rank of the module, which is
irrelevant to the topology. The actual obstruction group is the reduced group $\tilde K_0$,
defined to be either the cokernel of the inclusion of the trivial group, $K_0(Z[\{1\}])\to
K_0(Z[\pi])$, or the kernel of the rank homomorphism $K_0(Z[\pi])\to Z$. 

\subhead 2.4 Tameness\endsubhead 
Mapping cylinder neighborhoods have special homotopy properties. The eventual result is that
certain of these actually characterize mapping cylinders, modulo  $K$-theory problems. We
describe these.

The first property is that the neighborhood deformation retracts to $Y$, by pushing toward
the 0 end of the mapping cylinder. A key feature of this deformation is that the complement
of $Y$ stays in the complement until the last instant when everything collapses into $Y$. We
formalize this as: 
\subsubhead Definition\endsubsubhead An embedding $Y\subset X$ is  {\it forward tame\/} if
there is a  map $f\:X\times I\to X$ satisfying
\roster \item $f(x,t)=x$ if $t=0$ or $x\in Y$;
\item $f^{-1}(Y)=Y\times I\cup U\times\{1\}$, where $U$ is some neighborhood of $Y$.
\endroster
On the other hand we could pull the complement of $Y$ away from $Y$ by pushing toward the
other end of the mapping cylinder.  This formalizes to:
\subsubhead Definition\endsubsubhead An embedding $Y\subset X$ is  {\it backwards tame\/} 
if  there is  a map $b\:(X-Y)\times I\to X-Y$ satisfying
\roster\item $b(x,t)=x$ if $t=0$; and
\item for every $t>0$ the closure in $X$ of $b((X-Y)\times\{t\})$ is disjoint from  $Y$. 
\endroster
Putting these together we say:
\subsubhead Definition\endsubsubhead An embedding $Y\subset X$ is {\it tame\/} if it is both
forward and backward tame. 

Quite a bit is known about tameness. For instance if the embedding has finitely presented
constant local fundamental groups (see below) then there are homological characterizations,
and forward and backward tameness are equivalent because their homological formulations are
Poincar\'e dual
\cite{QS, 2.14}. If the embedding has {\it trivial\/} local fundamental groups then it is
always tame because the homological conditions are implied by the ANR hypotheses and
excision \cite{QS, 2.12}. See \cite{HR} for a treatment in the nonmanifold case.

\subhead 2.5 Homotopy links\endsubhead 
One of the main applications of tameness is to give a
comparison of the embedding with a ``universal'' mapping cylinder constructed using the
homotopy link. 
\subsubhead Definition\endsubsubhead The {\it homotopy link\/} of $Y\subset X$, denoted $\t
holink (X,Y)$, is a subset of the space of paths in $X$ with the compact-open topology.
Specifically it consists of the paths $s\:[0,1]\to X$ with $s^{-1}(Y)=\{0\}$. Evaluation at
0 gives a map $\t ev _0\:\t holink (X,Y)\to Y$. The whole evaluation map is $\t holink
(X,Y)\times I\to X$, and continuity implies this factors through  a map on the mapping
cylinder;
$\t ev \:\t cylinder (\t ev _0)\to X$. This preserves complements and is the identity on $Y$,
and in fact is the universal such map from a mapping cylinder. 

Now suppose $Y$ has a mapping cylinder neighborhood $N\simeq \t cylinder (q)$, with map
$q\:\partial N\to Y$. Since the homotopy link cylinder is universal, the geometric one
factors through it. More explicitly, each point in $\partial N$ determines a cylinder arc in
$N$. These arcs are points in the homotopy link so define a map $\partial N \to \t holink
(X,Y)$. This extends to a map of mapping cylinders. Further, $\partial N \to \t holink
(X,Y)$ turns out to be a controlled homotopy equivalence over $Y$, so the homotopy link
provides a homotopy model for any geometric mapping cylinder neighborhood.

Some of this last construction can be done using forward tameness in place of an actual
mapping cylinder. Suppose $f\:X\times I\to X$ is a forward-tameness deformation, and $U$ is
a neighborhood of $Y$ with $f(U\times\{1\})\subset Y$. Then the arcs $f\:\{x\}\times I\to X$
for $x\in U-Y$ define points in the homotopy link. This defines a map $U-Y\to \t holink
(X,Y)$. Using this in the first coordinate and distance from $Y$ in the second gives a map
to the universal mapping cylinder,
$U\to \t cylinder (\t ev _0)$. When $Y$ is also backwards tame this map is in an
appropriate sense a controlled local equivalence near $Y$. Tameness therefore encodes
essentially the same local homotopy information as a mapping cylinder neighborhood. 

 \subhead 2.6 Controlling local fundamental groups\endsubhead 
In standard (uncontrolled) geometric topology the fundamental group plays a central role.
Roughly this is because algebraic topology is effective with 1-connected spaces, and general
spaces are made 1-connected by taking universal covers. In controlled topology the same
principle applies, but fundamental groups cannot
be used directly because (among other problems) their definition depends on choices of
basepoints. Instead we use comparisons with reference spaces.

The general setting is a reference map $p\:E\to Y$, the controlled thing being studied, $W\to
Y$, and a  map $f\:W\to E$ that is required to commute with maps to $Y$ up to some
error $\delta$.  $f$ is said to be
$(\delta, 1)$-connected if given  a relative 2-complex $(K,L)$ and  a
$\delta$-commutative diagram 
$$\CD L@>>> W\\
@VV{\subset}V @VV{f}V\\
K@>>>E\endCD$$
then there is an extension $K\to W$ whose composition into $Y$ is within $\delta$ of $K\to
E\to Y$. When this is satisfied $W$ and $E$ have the same local fundamental group behavior
over $Y$, even if ``local fundamental groups'' do not make sense.

The importance of using reference maps to control $\pi_1$  increases with
increased complexity of local $\pi_1$ behavior. If the geometric situation is locally
1-connected over
$Y$ then no
$\pi_1$ control is needed. If the local fundamental groups are constant then a locally
1-connected covering space can be used. If local fundamental groups are locally constant over
$Y$ then we can use covering spaces of inverse images of open sets in $Y$. But now the
situation starts getting complicated: in geometric constructions we are controlling sizes,
so we need to know  these open sets are fairly large. In fact we need a priori estimates on
these sizes so geometric data can be chosen small in comparison. The simplest way to do
this is to control local $\pi_1$ using a fixed reference map $E\to Y$. In this way
whatever size data we need is determined by the reference map, and doesn't have to be made
explicit to be controlled. In the most general situation local fundamental groups change
from place to place. This is easily encoded using reference maps and awful to do with
groups. 

\subhead 2.7 Stratified systems of fibrations\endsubhead In the previous section we
described reference maps as a way to avoid the awkwardness of group formulations of local
$\pi_1$ structure.  However geometric constructions do use group formulations. Core steps of
proofs are usually done assuming constant local fundamental groups and using locally
1-connected covers.  General cases are obtained from this by fitting together locally
constant pieces. Thus the general
$\pi_1$ control apparatus is not intended to feed directly into core proofs, but to 
formulate general hypotheses that in proofs  inductively reduce to constant
cases. ``Stratified systems of fibrations'' \cite{QE2} work well for this.

\subsubhead Definition\endsubsubhead Suppose $p\:E\to Y$ is a map, and $Y=Y^n\supset
Y^{n-1}\supset\dots\supset Y^0$ is a filtration by closed subsets. $p$ is a {\it stratified
system of fibrations\/} (with filtration $Y^*$) if
\roster\item the restriction to each of the strata, 
$$p^{-1}(Y^i-Y^{i-1})@>p>> Y^i-Y^{i-1}$$
is a fibration, and
\item each term in the filtration is a $p$-NDR. This means there is a neighborhood of $Y^i$,
a deformation of it into $Y^i$ in $Y$ that preserves strata until the very end, and this
deformation is covered by a deformation of the inverse image in $E$. 
\endroster

Lots is known about these. There are many examples, reductions to apparently weaker data,
constructions, etc., see \cite{QS, CS, H}. 

In the mapping cylinder context the tameness hypothesis  provides us with a
canonical reference map, the homotopy link. Local $\pi_1$ hypotheses are formulated
in terms of this. Standard procedure (see the statement in 2.10) is to assume there is a 
stratified system of fibrations  $E\to Y$ and a map  $\t holink (X,	Y)\to E$ that is 
locally 1-connected over
$Y$. In many applications the homotopy link itself is a stratified system of fibrations.

 \subhead 2.8 Homology \endsubhead 
The mapping cylinder problem has obstructions lying in locally finite homology with
coefficients in a spectral cosheaf. This sounds complicated but is actually good news:
nothing simpler could work; it is reasonably accessible to calculation; and the formal
properties alone have important applications. In this section we outline the general setup
developed in \cite{QE2}. We assume general familiarity with the use of spectra to construct
homology theories. 

The basic setting is a spectrum-valued functor of maps with locally-compact target.
In more detail, the domain of this functor is the category with objects maps $p\:E\to B$, with
$B$ a locally compact metric space. Morphisms are pairs of maps $(F,f)$ forming a
commutative diagram 
$$\CD E_1@>F>>E_2\\
@VV{p_1 }V@VV{ p_2 }V\\
B_1 @>f>> B_2\endCD$$
and so that $f$ is proper (inverse images of compact sets are compact). In the application
$B$ is the control space where sizes are measured, and $E$ serves to control local
fundamental groups as in Section~2.5. 

We explain the local compactness hypothesis. The
technical work concerns manifolds mapping to $E$. We work over small open sets in $B$, and the
inverse image must have compact closure in the manifold. To get this we assume the map from
the manifold to $B$ is proper. But then we have to restrict to proper maps of $B$ to
preserve this property. We cannot simply require the manifold to be compact because we
need a restriction operation that destroys compactness. If $U\subset B$ is an open set then
restriction to inverse images of $U$ gives a map from manifold gadgets over $B$ to ones over
$U$. Even if we start with a compact manifold over $B$ the result will usually be noncompact
(but proper) over $U$. 

The homology of $X$ with coefficients in a spectrum $\b J $ is  the spectrum
$X\wedge \b J $. In more detail, $\b J $ is a seqence of spaces $J_n$ with various maps.
Start with the sequence of spaces $X\times J_n$, divide out $X$ times the basepoint in
$J_n$. The maps for $J_*$ then give this sequence of spaces the stucture of a suspension
spectrum. The associated spectrum is $X\wedge \b J $. We also denote this by $\b H (X;\b J
)$. Note this is a spectrum; the homology {\it groups\/} are the homotopy groups of this
spectrum
$$H_i(X;\b J ) = \pi_i(\b H (X; \b J ).$$
Note also that (unlike ordinary homology) these groups may be nontrivial for $i<0$. 

The locally-compact wrinkle in the theory requires us to work with locally finite homology.
This is essentially the relative homology of the 1-point compactification.

The ``Atiyah-Hirzebruch'' spectral sequence (due originally to G\. Whitehead) is a spectral
sequence of the form 
$$E^2_{i,j}=H_i(X;\pi_j\b J )\implies H_{i+j}(X; \b J ).$$
From this one sees, for example, that $H_j(X;\b J )$ always vanishes for $j<j_{\t min }$
exactly when
$\pi_j(\b J )=0$ for $j<j_{\t min }$, and that groups near the vanishing line are
quite accessible. This turns out to be very useful in applications. 

We now return to the context of a spectrum-valued functor $\b J (p)$, defined on the
category of maps $p$ with locally compact metric range spaces. In this case we can define a
``sheaf'' generalization of the homology construction. Suppose $p\:E\to Y$ is a map in the
category. We can apply the functor fiberwise to get a spectrum $\b J (p^{-1}(y)\to y)$ over
each $y\in Y$. With  mild additional assumptions we can fit these
together to get a ``spectral cosheaf'' over $Y$. This is a sequence of spaces
$J_n(p^{-1}(\#))$ with maps to $Y$ and maps to each other making the fibers over $Y$ into
spectra. In the constant-coefficient case $F\times Y\to Y$ this just gives $J_n(F\to \t pt
)\times Y$. By analogy with the constant-coefficient case we define {\it homology } with
coefficients in this cosheaf by first dividing out the 0-section of each $J_n(p^{-1}(\#))\to
Y$, then taking the spectrum associated to the resulting suspension spectrum. We use the
notation
$\b H (Y;\b J (p^{-1}(\#)))$ for this spectrum, and $H_*$ for its homotopy groups.

Again we actually need locally-finite homology. The spectrum for this is obtained by adding
a point at infinity to the spectral cosheaf, over the point at infinity in the 1-point
compactification of $Y$. Then divide out the 0-section and proceed as before. If $Y$ is
already compact this does not change the homology.

There is a generalization of the Atiyah-Hirzebruch spectral sequence to the non-constant
coefficient case. Namely in the situation of the previous paragraph we get 
$$E^2_{i,j}=H_i(X;\pi_j\b J (p^{-1}(\#)))\implies H_{i+j}(X; \b J (p^{-1}(\#))),$$
where the groups on the right are ``ordinary'' cosheaf homology groups. Again these are
reasonably accessible near the vanishing line for the coefficient spectra.

The is a  useful extension of the spectral cosheaf construction. Suppose,
as before, that
$E\to Y$ is a map in the domain of the functor, but now assume also that $f\:Y\to Z$ is a
proper map. Then we can construct a spectral cosheaf over $Z$ by applying the functor to
inverses under
$f$. More explicitly, over a point $z\in Z$ we put the spectrum $\b J (p^{-1}(f^{-1}(z))\to
f^{-1}(z))$. As before we can define a homology spectrum by dividing out 0-sections and
taking the associated spectrum. The output of this construction is denoted by $\b H
(Z;\b J (p^{-1}(f^{-1}(\#))))$. The notation is a bit tricky. Note we can do the previous
construction to the composition $fp$ and get a spectral cosheaf denoted $(fp)^{-1}(\#)$.
This is different from the cosheaf just constructed, though in some cases they have the same
homology.

\subhead 2.9 Assembly maps\endsubhead  We continue with the terminology of the
previous section. Suppose $p\:E\to Y$ is an object in the  category of ``proper maps
to locally compact spaces''. Then for each $y \in Y$ the inclusion
$$\CD p^{-1}(y)@>>> E\\
@VV{p}V @VV{p}V\\
y@>{\subset}>>Y\endCD$$
is a morphism in the category.  Applying $\b J $ gives maps from fibers of the spectral
cosheaf into $\b J (p)$.  Under mild continuity hypotheses these fit together to 
give a map on the total space of the cosheaf. Since the target of this map is a spectrum the
map factors through the associated spectrum of the total space to define a map of spectra
$$\b H (Y;\b J (p^{-1})) @>>> \b J (p).$$
This is the ``general nonsense'' description of the assembly map. In special cases there are
other descriptions that may give better understanding.

We will make use of the functoriality of assembly maps. Suppose there is a morphism 
$$\CD E_1@>F>>E_2\\
@VV{p_1}V@VV{p_2}V\\
Y_1@>f>>Y_2\endCD$$
in the domain category of
 $\b J $. Then  following diagram of spectra commutes:
$$\CD \b H (Y_1;\b J (p_1^{-1})) @>>> \b H (Y_2;\b J (p_2^{-1}))\\
@VVV@VVV\\
\b J (p_1)@>>>\b J (p_2)\endCD$$
The top map also factors through the mixed homology $ \b H (Y_2;\b J (f(p_2^{-1}))$.

\subhead 2.10 $K$-theory\endsubhead 
The previous section gives the context for homological obstructions. In this section we
discuss particular functors used to make contact with the topological problems. 

The logical context for the next theorem is that geometric-topological techniques can be used
to formulate obstruction groups for controlled problems (see Section 3). This tells us what
they are good for, but says very little about their nature. The next theorem provides another
description that displays  global properties. This is incorporated in the final
statement in~\S3.1.
\proclaim{Controlled assembly isomorphism theorem} There is a spectrum-valued functor $\Cal S$
defined on maps to locally compact metric spaces, such that if $p\:E\to Y$ is a stratified
system of fibrations over a locally compact finite-dimensional ANR then
\roster \item $\pi_0\Cal S(p)$ is the obstruction group for mapping cylinder neighborhoods
of $Y$ with local fundamental groups modeled on $p$, and
\item the assembly map $\b H ^{lf}(Y;\Cal S(p^{-1}(\#))\to \Cal S(p)$ is an equivalence of
spectra.
\endroster
\endproclaim
First we explore the significance of conclusion (1). Uncontrolled work determines
some of the homotopy of the ``coefficient'' spectra ($Y$  a point):
$$\pi_i(\Cal S(F\to pt))=\cases Wh(\pi_1F)& \text{\quad when }i=1;\\
\tilde K_0(Z[\pi_1F])& \text{\quad when }i=0,\text{ and}\\
K_{-i}(Z[\pi_1F])& \text{\quad when }-i<0.\endcases
$$
The $i=1$ case comes requiring the same spectrum to work for h-cobordisms, $i=0$ is the
uncontrolled end theorem (Siebenmann, see \S2.3), and $-i<0$ comes from seeing Bass'
definition of lower
$K$-theory \cite{B} come out of tinkering with Euclidean spaces $Y=R^i$ (\cite{PW}). One can
also require a connection to pseudoisotopy and get $\pi_2$ to be $Wh_2$, \cite{QE4}. 
The spectral sequence shows the higher homotopy plays no role in the
obstructions of interest, so we don't
particularly care what it is. 

There are many constructions of spectra encoding lower $K$-theory, and many of these extend
to spectrum-valued functors satisfying condition (1) of the theorem. Conclusion (2), which
is the source of the real power of the theory, is much more delicate. 

The proof of (2) follows  the proof of uniqueness of homology,  i.e\.  a morphism
of homology theories that induces an isomorphism on homology of a point  is an
isomorphism on finite-dimensional ANRs. The proof proceeds by induction using exact
sequences, first establishing isomorphism for spheres, then finite CW complexes, then (by a
trick) ANRs. To prove (2) this way we (i) show the right side ($\Cal S(p)$) satisfies
appropriate versions of the axioms of homology in the $Y$ variable; (ii) observe that the
map  gives isomorphisms over points by definition; and (iii) make minor adjustments to
incorporate the ``coefficient system'' (reference map $p$). The hard part of this is (i).
The axioms are homotopy; excision; and a fibration condition for pairs. The fibration
hypothesis is the spectrum version of the long exact sequence of homology {\it groups\/} of
a pair: the long exact sequence is the homotopy sequence of the fibration. Technically since
we are working with locally finite homology the pair axiom is replaced by a condition on
restrictions to open sets, but it amounts to the same thing. 

Several conclusions can be drawn from this outline. First, we may not care about the higher
homotopy of the spectrum but Nature does. The inductive plan of the proof only works if
{\it all\/} the groups line up correctly, so we have to get them right whether we want to use
them or not. Second, the proof is all-or-nothing. Again because it depends on an induction
 it either works or fails, and there are no interesting partial results when it
fails. Finally, the bottleneck in such arguments is usually the excision axiom. Recall that
this is the one that fails for homotopy groups, and thus enables homotopy to be so much more
complicated than homology.

The theorem is proved in \cite{QE2} with a spectrum $\Cal S$
constructed  using pseudoisotopy of manifolds. This is a version of what is now known as
$A$-theory or Waldhausen $K$-theory. For current topological applications one such spectrum
is enough. However the proof in \cite{QE2} is complicated and not too clear, so there have
been efforts to find other proofs of this key step. Also, a version extending algebraic
$K$-theory would have significant applications to algebra. So far this has been not been
done: none of the other formulations of Waldhausen $K$-theory and none of the algebraic $K$
constructions (Quillen, Volodin and so on) have been acceptable to the methods of the proof. 
The author thinks he  has  a construction for algebraic $K$-theory, but he has 
thought this before (cf\.
\cite{QK}) so skepticism is appropriate until details appear.

We offer a philosophical explanation of why the controlled assembly isomorphism theorem is
so hard to prove. Frequently a complicated proof is ``explained'' by the existence of a
false similar statement. For example the ``reason'' Freedman's topological embedding theorem
for 4-manifolds is so hard is that the analogous statement for smooth embeddings is false.
The proof must have a topological construction so bizarre that it cannot possibly give a
smooth outcome, and must depend on it so essentially that it cannot possibly be avoided.
What then is the false thing forcing the isomorphism theorem proof to be so intolerant?
 The problem is probably in {\it quadratic\/}
stable algebra (surgery, $L$-theory). There the geometrically significant lower homotopy
groups of the spectrum do have some imprint of strange behavior in the higher groups. As a
result the surgery spectrum constructions now known {\it cannot\/} satisfy the isomorphism
theorem. The proof must be delicate enough to reject these impostors. Apparently $K$-theory
doesn't do anything strange enough to deserve the complexity; it is just an innocent victum
of problems in surgery. 
\head 3. The theorem\endhead 

Collecting the hypotheses developed in \S2, we suppose $X$ is a locally
compact finite dimensional metric ANR,
$Y\subset X$ is tame, $X-Y$ is a
manifold, $p\:E\to Y$ is a stratified system of fibrations, and there is
a controlled 1-connected map $\t holink (X,Y)\to E$.

\proclaim{3.1 Mapping cylinder existence theorem}  Under these conditions there is an
invariant 
$q_0(X,Y)\in H^{lf}_0(Y;\Cal S(p^{-1}(\#)))$. This vanishes if\/ $Y$ has a mapping
cylinder neighborhood, and conversely if the invariant vanishes and\/  $\t dim \,X-Y\geq 6$
then there is a mapping cylinder neighborhood.
\endproclaim
\subsubhead Dimension 5\endsubsubhead
This is still true in dimension  5  when the local fundamental groups of
$p$  have subexponential growth \cite{FQ, KQ}.

We outline the proof only well enough to show the major features.
\subhead 3.2 Nice neighborhoods and the obstruction \endsubhead
The key objective is to find neighborhoods $N$ with the right controlled homotopy type:
closed manifold neighborhoods so that $\partial N\to N-Y$ is an $\epsilon$ homotopy
equivalence over $Y$. Tameness is the main ingredient. Choose any small manifold
neighborhood $N$, and choose handlebody structures. The tameness deformations provide
homotopy data to show how to swap handles to make $\partial N\to N-Y$ highly connected. The
final step, which would make it a homotopy equivalence, is obstructed. In the uncontrolled
case we see a single nonvanishing relative homology group. If it is stably free over the
group ring then we can stabilize and swap handles corresponding to a basis to get a good $N$.
This module is a direct summand of a finitely generated free chain group, so is finitely
generated projective. The obstruction is the equivalence class of the projective module,
modulo stably free modules. In other words, its image in $\tilde K_0(Z[\pi_1])$. We indicate
modifications needed in the controlled setting. We can't use homology because this is a
quotient and quotients destroy size estimates. Instead we directly use the projection on
controlled chain groups. We define $\tilde K_0(Y;p,\epsilon,\delta)$ to be free modules over
$Y$ with ``$Z[\pi_1(p^{-1}(\#))]$'' coefficients (this is clarified in \S3), with projections
of radius $<\delta$, modulo ones with basis-preserving projections of radius $<\epsilon$.
Adding estimates to the uncontrolled argument gives an element of this set, and shows that
if it is trivial then the argument can be completed to get a nice $N$. 

We pause the proof to expand on the obstructons. 
The main point is  that although we have an ``obstruction'', and can arrange
for the set 
$\tilde K_0(Y;p,\epsilon,\delta)$ in which it lies to be a group, we know nothing about it.
This is where the characterization theorem of \S9 takes over. This shows:
\roster\item these groups are stable in the sense that for every $\epsilon>0$ there is
$\delta>0$ so that the map from the inverse limit 
$$\tilde K_0(Y;p,\epsilon,\delta) @<<< \lim_{\leftarrow}\tilde K_0(Y;p,*,*)$$
is an isomorphism; and 
\item the inverse limit is the spectral sheaf homology group $H^{lf}_0(Y;\Cal S(p^{-1}(\#)))$.
\endroster
The stability in (1) is  subtle and actually harder to prove than the description of
the limit in (2). For instance in the controlled algebra used here, sizes grow when
morphisms are composed. This means morphisms of fixed size do not form a category, and in
place of the homological and categorical techniques of the uncontrolled theory we have to
work with chain complexes and constantly estimate sizes. In contrast it is
possible to set up the inverse limit theory directly so the work takes place in a category
\cite{P}. The setup is more elaborate, but no estimates are needed and the group
$\lim_{\leftarrow}\tilde K_0(Y;p,*,*)$ appears as ordinary 
$K$-theory of a category. In some applications the stability property is essential (see
\S5). However for mapping cylinders it is not. We have extracted an invariant from a single
sufficiently small neighborhood $N$. But one could repeat the construction at smaller scales
to get a sequence of neighborhoods $N_i$ with estimates going to 0. From this we could
extract a sequence of related algebraic objects with estimates going to 0, or in other words
an element of the inverse limit. This approach may yet have significant applications.
However so far the benefits (convenience for the categorically sophisticated) do not seem to
outweigh the drawbacks (weaker theorems, more elaborate setups). 

\subhead 3.3 Getting mapping cylinders\endsubhead
Returning to the proof, we suppose the obstruction vanishes so we can find nice
neighborhoods $N$. Repeat at smaller scales to get a decreasing sequence $N_i\supset
N_{i+1}\dots$ which are ``nice'' with decreasing size estimates. Recall that ``nice'' meant
roughly that the inclusion $\partial N_i\to N_i-Y$ is a controlled homotopy equivalence.
It follows that the regions between these are controlled h-cobordisms. Explicitly, the
inclusions of $\partial N_i$ and $\partial N_{i+1}$ in $N_i-\t interior (N_{i+1})$ are
controlled homotopy equivalences. If these h-cobordisms are all products then we can fit
together product structures $\partial N_i\times[\frac1{i+1},\frac1i]\simeq N_i-\t interior
(N_{i+1})$ to get a product structure $\partial N_1\times (0,1]\simeq N_1$. The control on
the size of the product structures shows the images of the arcs converge in $Y$, so this
gives a mapping cylinder. 

The intermediate regions originally constructed may not  be products, but we can
use a ``swindle'' to make them so. If we factor each $N_i-\t interior (N_{i+1})$ as a
composition of h-cobordisms $U_i\cup V_i$, $N_1$ becomes an infinite union $(U_1\cup
V_1)\cup (U_2\cup V_2)\cup \dots$. Reassociating expresses it as $U_1\cup
(V_1\cup U_2)\cup\dots$. The idea is choose the decompositions so the new pieces, $V_i\cup
U_{i+1}$ are all products. In the uncontrolled setting this is a simple consequence of the
invertibility of h-cobordisms. The controlled version is not so simple. We want to
inductively choose the decomposition
$U_{i+1}\cup V_{i+1}$ so the union $V_i\cup U_{i+1}$ is a product. But we must maintain
finer control on $U_{i+1}$ than is available on $V_i$. The argument thus uses stability of
h-cobordism obstructions: we need not only that $V_i$ has some inverse, but that it has one
with arbitrarily finer control. This is a deep fact, so this ``swindle'' is not just a
formal argument. As explained above this can be avoided by working  with a sequence
$\{N_i\}$ to formulate the obstruction directly as an element of the inverse limit. When
this vanishes the intermediate regions are automatically already products. 

This completes the sketch of the proof. 
\subhead 3.4 Refinements\endsubhead
We give two refinements that follow from the proof. The first concerns smooth or PL
structures, and the second  provides a way to recognize
mapping cylinders themselves, not just existence of neighborhoods.

\proclaim{Smooth and PL cylinders} If the manifold in Theorem 3.1 has a smooth or PL structure
then the mapping cylinder can be chosen to be smooth or PL, in the sense that the submanifold
$N-Y$ is, and the map
$\partial N\times (0,1]\to N-Y$ is a diffeomorphism or PL isomorphism.
\endproclaim

The proof of the Theorem uses handlebody
theory, which works in any category of manifolds (dim
$>4$ in the topological case). Thus the argument and obstructions are category-independent: if
$X-Y$ has a smooth structure we get a smooth $N$, etc. 

There is also a  smoothing and triangulation theory  that shows a topological
mapping cylinder in a PL manifold can be made PL, and similarly for smoothing. Using this we
could deduce the structure refinement from the topological case. The point of observing  it
 directly from the proof
 is that eventually it is possible to run the argument
backwards and {\it derive\/} the smoothing and triangulation structure theory from 
controlled theorems. In such ways the controlled theory unifies as well as extends the older
work.

The structure refinement is {\it not\/} true in dimension 5, no matter how nice the local
fundamental groups are.

The second result of the section gives a criterion for $X$ itself to be a mapping cylinder
over $Y$. As in the discussion of tameness in 2.4 we extract properties of the radial
deformation of a mapping cylinder. If $X=\t cyl (g)$ for some map $N\to Y$ then the radial
deformation to $Y$ is a map $f\:X\times I\to X$ satisfying:
\roster\item $f^{-1}(Y)=Y\times I \cup X\times \{1\}$
\item $f(x,t)=x$ if $t=0$ or $x\in Y$, and
\item $f(f(x,t),1)=f(x,1)$
\endroster
The last condition means that if we use the time-1 retraction $f_1\:X\to Y$ as a control
map then the deformation $f$ has radius 0 in
$Y$. The criterion relaxes this, requiring only that $f$ has radius less than some
appropriate $\delta$. Note that to be useful this ``appropriate $\delta$'' must be known in
advance, before $X$ and $f$ are chosen. Note also that if $Y$ is not compact then this sort
of control uses a {\it function\/} $\delta\:Y\to (0,\infty)$. Typically these functions go
to 0, so provide progressively finer control near the ends of $Y$.

\proclaim{Mapping cylinder recognition} Suppose $Y$ is a locally compact finite dimensional
ANR, 
$p\:E\to Y$ is a stratified system of fibrations, and a dimension $n\geq 6$ is given. Then
there is $\delta>0$ so that if
\roster\item $X\supset Y$ with $X-Y$ a manifold (with boundary) of dimension $n$;
\item there is a map $\t holink (X,Y)\to E$ that is $(\delta,1)$-connected over $Y$;
\item $f\:X\times I\to X$ is a proper deformation retraction of $X$ to $Y$ that preserves the
complement of $Y$ when $t<1$, and $f_1f$ has radius $<\delta$ in $Y$.
\item the inclusion $\partial X\subset X-Y$ is $(\delta,1)$-connected over $Y$, using $f_1$
as control map.
\endroster
Then there is a map $g\:\partial X\to Y$ and a homeomorphism $\t cyl (g)\to X$ that is
 the identity on the boundary and $Y$.
\endproclaim
We can further arrange for the radial deformation in $\t cyl (g)$ to be close to the
deformation $f$. More specifically choose $\epsilon>0$ along with $Y$, $p$, and $n$. Then
there is a choice of $\delta$ so we can get $g$ and  homeomorphism $h\:\t cyl (g)\to X$
with the diagram 
$$\CD \t cyl (g)\times I@>{h\times \t id }>> X\times I\\
@VV{\t {cyl. projection} }V@VV{f_1f}V\\
Y@>\t id >>Y\endCD$$
commutative within $\epsilon$. 

The hypotheses collected in this theorem encode the properties of the ``nice neighborhoods''
used in the proof of the theorem of 3.1. The part of the proof outlined in 3.2  shows the
obstruction vanishes if and only if nice neighborhoods exist. The argument in 3.3 then
proves the theorem stated here, that nice neighborhoods are mapping cylinders. 

We remark on the role of stability. The argument in 3.3 requires finding a descending
sequence of nice neighborhoods with size parameters going to 0. These exist because the
obstructions are stable (the same at all sufficiently small scales), and the initial nice
neighborhood is chosen to have size in the stable range so it's existence shows the
obstructions are trivial. The abstract existence of mapping cylinder neighborhoods can be
formulated directly in terms of the inverse limit, avoiding stability. However the
recognition theorem is not accessible to this approach because the existence of a single nice
neighborhood does not show the vanishing of the obstruction.

\head 4. Applications\endhead 
The applications detailed here are briefly described in the introduction.

\proclaim{4.1 Mapping cylinders in manifolds \cite{QE1}} Suppose $Y\subset M$ is a closed
embedding with locally 1-connected complement, of an ANR in the interior of a manifold of
dimension
$\geq5$. Then
$Y$ has a mapping cylinder neighborhood.\endproclaim
To derive this from the main theorem we must show tameness and vanishing of the
obstruction. The homological characterization of forward tameness
follows from excision and triviality of local fundamental groups. Backwards tameness follows
from this. Since the local fundamental groups are trivial we can use the identity
$Y\to Y$ as the control map. The main theorem identifies the obstruction as lying in
$H^{lf}_0(Y;\Cal S(id))$. The spectral sequence for this has $E^2$ terms $H^{lf}_0(Y;\tilde
K_0(Z))$ and   $H^{lf}_i(Y;
K_-i(Z))$ for $i>0$. But $\tilde
K_0$ and the lower $K$-theory of $Z$ is all trivial, so the obstruction group is trivial.
Therefore a mapping cylinder exists.

\proclaim{4.2 Finiteness of ANRs} A compact finite-dimensional ANR is homotopy equivalent to
a finite complex.\endproclaim
This statement is the ``Borsuk conjecture'', and was proved for all compact ANRs (not
necessarily finite-dimensional) by J. West. The finite-dimensional case follows from the
previous result as follows: a compact finite-dimensional ANR has an embedding in some
Euclidean space. If this is has locally 1-connected complement then there is a mapping
cylinder neighborhood. The neighborhood is a manifold (smooth, actually) so is a finite
complex. The mapping cylinder projection is a homotopy equivalence. What if the embedding
does not have locally 1-connected complement? The inclusion into $R^{n+1}$ has locally
0-connected complement, and if an embedding has locally 0-connected complement then the
inclusion into $R^{n+1}$ has locally 1-connected complement. Thus increasing dimension by
two always makes the complement locally 1-connected.

\proclaim{4.3 Collaring in homology manifolds} Suppose $M$ is an ANR homology manifold of
dimension $\geq 5$, the embedding $\partial M\subset M$ has locally 1-connected complement,
and
$M-\partial M$ is a manifold. Then there is a mapping cylinder neighborhood of $\partial M$.
There is a collar (neighborhood homeomorphic to
$(\partial M)\times I$) if and only if $(\partial
M)\times R$ is a manifold.\endproclaim
The inclusion of the boundary of a homology manifold is homologically locally
infinitely-connected, but local fundamental groups may be nontrivial. For instance
closure of the strange component of the complement of the Alexander horned sphere in $S^3$ is
a homology manifold with non-locally 1-connected complement of the boundary. With the
1-connected hypothesis the proof of mapping cylinders is the same as the previous theorem
with a little modification of the proof of forward tameness. 

We give some context for the collaring statement. First, the map in the mapping cylinder
must be a resolution (map from a manifold to a homology manifold that is a local homotopy
equivalence, see \S4). Edwards' resolution theorem is that when the dimension is $\geq5$ this
map can be approximated by a homeomorphism if and only if the homology manifold has the
``disjoint 2-disk property.'' If $N$ is a homology manifold then an easy argument shows
$N\times R^2$ has the disjoint 2-disk property so $N$ resolvable implies $N\times R^2$ is
resolvable, and therefore a manifold. It is one of the outstanding conjectures in the area
that $N\times R$ is already a manifold. The theorem shows this conjecture is equivalent to
existence of collars of boundaries in certain homology manifolds. 

The proof of the collaring statement follows from the  fact that the mapping cylinder map is
a resolution, so Edwards' theorem shows the product with $R$ can be approximated by a
homeomorphism if and only if $(partial
M)\times I$ is a manifold.

\subhead 4.4 Stratified spaces\endsubhead Many interesting spaces are
not manifolds, but are built of manifold pieces. These include algebraic varieties,
stratifications coming from singularities, and polyhedra. Generally a ``stratified
space'' 
 has a closed filtration
$X=X_n\supset X_{n-1}\dots\supset X_0$ and the strata $X_i-X_{i-1}$ are  required to
be manifolds. There are several versions that differ in the way the strata fit together.  The
geometric versions (Whitney, Thom, and PL)  have mapping cylinder neighborhoods and
complicated relations among them as part of their structure. The most successful topological
version, homotopy, or ``Quinn'' stratified spaces \cite{QS}, were identified as an outgrowth
of controlled topology and have local homotopy conditions  relating the strata. In these the
strata may not have mapping cylinder neighborhoods, and the obstruction is exactly the one
identified in 2.10. 

More specifically, suppose $X$ is a homotopy stratified space in the sense of
\cite{QS}. Then more-or-less by definition
\item{$\bullet$} the embedding $X_{i-1}\subset X_i$ is tame;
\item{$\bullet$} the projection of the homotopy link $ev_0\:\t holink (X_i,X_{i-1})\to
X_{i-1}$ is a stratified system of fibrations; and
\item{$\bullet$} $X_i-X_{i-1}$ is a manifold.

\noindent Thus we conclude there is an obstruction in $H^{lf}_0(X_{i-1};\Cal
S(ev_0^{-1}(\#)))$ to the existence of a mapping cylinder neighborhood of $X_{i-1}$ in
$X_i$. Vanishing of the obstruction implies existence of such a neighborhood if either $\t dim
X_i\geq6$ or
$\t dim X_i=4$ and fundamental groups of point-inverses in the homotopy link are ``good.''

We note that the fact that existence is obstructed means that mapping cylinders are not the
natural local structure in these spaces. A weaker version developed by Hughes and others seems
to be correct, see \S3.5.

\subhead 4.5 Topological  actions of finite groups\endsubhead Suppose a
finite group $G$ acts on a manifold $M$. We can filter the quotient $M/G$ by orbit types:
images of points lie in the same stratum of the quotient if their isotropy subgroups are
conjugate. If the action is smooth or PL then the quotient is a smooth or PL stratified
space, though the stratification may not be exactly the orbit type stratification. If the
action is just topological then really awful point-set things can happen in the quotient. A
nice compromise is the class of ``homotopically stratified'' actions \cite{QS}, where the
quotient with orbit type filtration is assumed to be homotopically stratified in the sense
discussed above. This rules out weird point-set behavior but allows many other things. For
instance these can have mapping cylinder problems. 

An ``equivariant mapping cylinder neighborhood'' of a $G$-invariant subset of $M$ is just
what it sounds like: a mapping cylinder structure invariant under the action of $G$. The
quotient is an ordinary mapping cylinder neighborhood of the quotient subset. Therefore
Theorem 3.1 can be applied in the quotient to determine the existence of equivariant mapping
cylinder neighborhoods. We discuss the easiest case, neighborhoods of the non-free points.

\proclaim{Theorem} Suppose the finite group $G$ acts in a homotopically-stratified way on a
compact manifold,  let $Y\subset M$ be the points not moved freely by $G$, and suppose the
codimension of
$Y$ is $\geq3$. Then 
\roster\item there is a stratified system of fibrations $p\:B_{G_{\#}}\to Y/G$ whose fiber
over $x\in Y$ is the classifying space of the subgroup $G_x\subset G$ fixing $x$;
\item there is an obstruction in $H_0(Y/G;\Cal S(p^{-1}(\#)))$ to the existence of an
equivariant mapping cylinder neighborhood of $Y$; and
\item if there is an equivariant cylinder neighborhood of $Y\cap\partial M$ in $\partial M$
and
$\t dim M\geq 5$, then the obstruction vanishes if and only if there is an extension of the
boundary structure to an equivariant mapping cylinder neighborhood of all of $Y$.
\endroster
\endproclaim
The hypothesis that $Y$ has codimension is at least 3 implies the embedding $Y\subset M$
has locally 1-connected complement. Therefore local fundamental groups in the quotient come
from the the group action, and are modeled by the isotropy groups described in the theorem.
The obstructions are often quite accessible:
\roster\item $K_{-i}(Z[H])=0$ for finite groups $H$ and $-i\leq -2$ \cite{C}. Therefore the
spectral sequence for the obstruction group has $E^2$ terms only $H_0(Y/G;\tilde
K_0(Z[G_{\#}]))$ and $H_1(Y/G;
K_{-1}(Z[G_{\#}]))$. (these are group, not spectral, cosheaf homology groups);
\item if $Y/G$ is connected and there is a point fixed by $G$ then the $H_0$ term reduces to
$\tilde K_0(Z[G])$; and
\item there is an action of a finite group on a disk that is smooth on the boundary and
locally linear (in fact can be smoothed in the complement of any fixed point), but
the non-free set does not have an equivariant mapping cylinder neighborhood because the
$\tilde K_0(Z[G])$ part of the obstruction is nontrivial
\cite{QE2}.
\endroster
\subhead 4.6 Topological regular neighborhoods\endsubhead 
Mapping cylinders are wonderful, but since they do not always exist they are not 
satisfactory objects for a topological theory of ``regular neighborhoods.'' 
The appropriate notion seems to be a skewed mapping cylinder neighborhood in
$X\times[0,\infty)$. 
\subsubhead Definition\endsubsubhead Suppose $X$ is locally compact and $Y\subset X$ is
closed. A {\it topological regular neighborhood\/} of $Y$ consists of
\roster\item an open neighborhood $U$ of $Y$,
\item a proper map $q\:U\to Y\times[0,1)$ that is the identity $Y\to Y\times\{0\}$ and
preserves complements of these sets,
\item a homeomorphism of the relative mapping cylinder $\t cyl (q,\t id _Y)$ with a
neighborhood of $Y\times \{0\}$ in $Y\times[0,\infty)$ that is the identity on
$Y\times[0,1)$ and $U\times\{0\}$.
\endroster
The relative mapping cylinder is the ordinary mapping cylinder with the cylinder arcs in the
subset $Y\times I$ identified to points. The result contains a copy of $Y$, and the
complement of this is the mapping cylinder of the restriction of $q$ to $U-Y\to Y\times
(0,1)$. The homeomorphism in (3) takes the cyinder arcs to arcs that start on $U\times\{0\}$
and go diagonally to  $Y\times[0,1)$, see the figure.

\figure1{Topological Regular Neighborhood}
The idea is that we may not be able to find  mapping cylinder neighborhoods because
there is an obstruction to finding an appropriate domain for the map. So we use the
neighborhood itself as the domain for a map, and get a mapping cylinder in the next higher
dimension. 

\proclaim{Existence of regular neighborhoods} Suppose $X$ is a locally compact ANR,
$Y\subset X$ is tame, and there is a map, controlled 1-connected over $Y$, from the homotopy
link of $Y$ in $X$ to a stratified system of fibrations over $Y$. Finally suppose $X-Y$ is a
manifold of dimension $\geq5$. Then there is a topological regular neighborhood of $Y$ in
$X$. \endproclaim

As usual this also holds for $X$ of dimension 4, provided the local fundamental groups have
subexponential growth.

Before indicating the proof we discuss some of the structure of these neighborhoods. Suppose
$B\subset A$ has a mapping cylinder neighborhood with map $q\:U\to B$. The homotopy link  is
in a sense universal for mapping cylinders mapping to $(A,B)$, so the cylinder
structure defines a map $U\to \t holink (A,B)$. This is an ``approximate fiber homotopy
equivalence'' over $Y$ \cite{QS, 2.7}. In the regular neighborhood situation the homotopy
link of $Y\times(0,1)\subset X\times(0,1)$ is the pullback of the homotopy link of $Y\subset
X$. Thus the regular neighborhood structure gives an approximate fiber homotopy equivalence
over $Y\times(0,1)$, $U\to \t holink (X,Y)\times(0,1)$. 

In the important special
case where $X$ is a homotopically stratified space as in \S4.4 and $Y$ is the union of the
lower strata, these regular neighborhoods are the same as the ``approximate tubular
neighborhoods'' developed by Hughes \cite{H}, and  earlier in special cases by  Hughes,
Taylor, Weinberger and Williams \cite{HTWW}.  In a stratified set $\t holink (X,Y)\to Y$ is
itself a stratified system of fibrations. Thus
$q\:U-Y\to Y\times(0,1)$ is approximately fiber homotopy equivalent to a stratified system
of fibrations. This identifies $q$ as a ``manifold stratified approximate fibration'' over
$Y\times(0,1)$.

Finally we outline how the theorem follows from a version of the mapping cylinder
recognition theorem of 3.4. Let $h\:X\times[0,1]\to X$ be a forward-tame deformation, and
let $V$ be a neighborhood of $Y$ such that $h(V\times\{1\})\subset Y$. Define
$f\:X\times[0,\infty)\times[0,1]\to X\times[0,\infty)$ by
$$f(x,s,t) =(h(x,t),s+t d(x,Y))$$
where $d(x,Y)$ is the distance from the point to the subspace $Y$. Then properties of $h$
imply
\roster\item $f$ is a homotopy from the identity at $t=0$ to a map at $t=1$ that takes
$V\times[0,\infty)$ into $Y\times[0,\infty)$;
\item when $t<1$ $f(\#,t)$ takes the complement of $Y\times[0,\infty)$ into itself;
\item $f(\#,t)$ is the identity on $Y\times[0,\infty)$;
\item $f(\#,t)$ takes the complement of $Y\times\{0\}$ into itself for all $t$. 
\endroster
Delete $Y\times\{0\}$, then this is a homotopy of $(X-Y)\times\{0\}\cup X\times(0,\infty)$.
We would like to arrange it to satisfy the conditions of the mapping cylinder recognition
theorem, to get a mapping cylinder over $Y\times(0,\infty)$. This cannot be done completely:
the end near $Y\times\{0\}$ is ok, but none of the conditions hold near $\infty$. Instead we
use a relative version: if the conditions hold over $Y\times(0,1)$ then some
neighborhood of
$Y\times\{0\}$ is a mapping cylinder over $Y\times(0,1-\epsilon)$. Such relative versions are
standard parts of controlled theory (see e.g\. the remarks before Theorem 1.3 in \cite{Q2}).
They are often not stated explicitly because the statements are so complicated, but
follow from the proofs. The actual goal is thus to arrange the conditions of the recognition
theorem to hold over $Y\times(0,1)$. 

Recall that we want to use $f_1$ as the control map. The first problem is that $f_1$ does
not even map all of the space into $Y\times(0,\infty)$. However this does work near 0.
Suppose the neighborhood
$V$ taken into $Y$ by $h$ contains the points within $\epsilon$ of $Y$. Then $f$ does deform
all of 
$f_1^{-1}(Y\times[0,\epsilon))$  into
$Y\times(0,\infty)$. Restrict to this, in the sense that we consider the space
$f_1^{-1}(Y\times[0,\epsilon))$ with control map $f_1$ over $Y\times[0,\epsilon)$.
Reparameterize $[0,\epsilon)$ as $[0,\infty)$. The situation is now that $f_1$ can serve as a
control map. It is also proper. We have lost something: since the original deformation did not
preserve
$f_1^{-1}(Y\times[0,\epsilon)$, the restriction does not define a deformation of the space.
However since $Y\times\{0\}$ is left fixed, by continuity
$f$ keeps some
$f_1^{-1}(Y\times[0,\tau)$ inside the new space. Reparameterize the interval again to arrange
the deformation to be defined on $f_1^{-1}(Y\times[0,2))$. We now have the control map and
deformation defined over $Y\times[0,2)$

The last step is to arrange  arbitrarily good size control, at least over $Y\times(0,1)$. 
The deformation $f$ is the identity on $Y\times\{0\}$, or in other words the composition
$f_1f$ has radius 0 as a homotopy of $Y\times\{0\}$ in itself. It follows that $f_1f$ has
very small radius over $Y\times[0,\epsilon)$, for small $\epsilon$. By reparameterizing
$[0,\epsilon)$ to $[0,2)$ we can arrange that $f_1f$ has arbitrarily small radius in the $Y$
coordinate over $Y\times[0,2)$. We need to do a little better. We are controlling over the 0
end of $Y\times(0,\infty)$, which is non-compact even if $Y$ is compact. Assume $Y$ is
compact to simplify the argument, then the control objective is a continuous
 function $\delta\:(0,\infty)\to (0\infty)$, not a constant $\delta>0$. Elaborate the
previous argument: since $f_1f$ has radius 0 over $Y\times\{0\}$, there is a continuous
increasing function
$\epsilon\:[0,\infty)\to[0,\infty)$ taking 0 to 0, so that $f_1f$ has radius $<\epsilon$ in
the $Y$ coordinate, over $Y\times[0,2)$. Now reparameterize by a homeomorphism
$\theta\:[0,\infty)\to[0,\infty)$ so that $\epsilon\theta <\delta$. The result is $\delta$
controlled in the $Y$ coordinate.  It remains to get control in the $[0,\infty)$ coordinate.
This is again a standard argument using continuity and reparameterization.

The outcome of all this is a neighborhood of $Y\times\{0\}$ in $X\times[0,\infty)$, a
control map to $Y\times[0,\infty)$, and a deformation defined over $Y\times[0,2)$ and
satisfying the control needed for the Recognition Theorem over $Y\times[0,1)$. The theorem
then asserts that there is a mapping cylinder structure provided the dimension of
$X\times[0,\infty)$ is at least 6, or in other words if $X$ has dimension at least 5, or $X$
has dimension 4 and the local fundamental groups are small.

\refstyle{A}

\bye